\documentclass
[12pt,a4paper]{article}

\usepackage{amssymb}
\usepackage{amsmath}

\newenvironment{dem}{\noindent {\bf Proof.}}{\hfill $\Box$\par}
\newtheorem{teorema}{Theorem}
\newtheorem{proposicion}[teorema]{Proposition}
\newtheorem{lema}[teorema]{Lemma}
\newtheorem{corolario}[teorema]{Corollary}

\newtheorem{ejem}[teorema]{Example}

\newtheorem{algoritmo}[teorema]{Algorithm}

\usepackage{color}
\newcommand{\WS}[1]{\textcolor{blue}{#1}}
\definecolor{verde}{rgb}{0,0.7,0}
\newcommand{\AS}[1]{\textcolor{verde}{#1}}

\newcommand{\bmx}{\left[\begin{matrix}}
\newcommand{\emx}{\end{matrix}\right]}

\begin{document}
\sloppy

\date{30 November 2010}

\title{
Feedback classification of linear systems over \mbox{von Neumann} regular rings}
\author{Andr\'es S\'aez-Schwedt and Wiland Schmale}
\maketitle

\begin{abstract}
It is proved that feedback classification of a linear system over
a commutative von Neumann regular ring $R$ can be reduced to the classification of a
finite family of systems, each of which is properly split into a reachable and a
non-reachable part, where the reachable part is in a Brunovski-type canonical form,
while the non-reachable part can only be altered by similarity.
If a canonical form is known for similarity of matrices over
$R$, then it can be used to construct a canonical form for feedback equivalence. An
explicit algorithm is given to obtain the canonical form in a computable context
together with an example over a finite ring.

\end{abstract}

\section{Introduction and motivation}

Let $R$ be a commutative ring with 1.
An $m$-input, {\it $n$-dimensional linear system over $R$} is a pair of matrices
$(A,B)$, with $A\in R^{n\times n}, B \in R^{n\times m}$.
The control theoretic background for this terminology is nicely
described in \cite[Introduction]{BBvV}. Within this context
two systems $(A,B)$ and $(A',B')$ are called {\it feedback equivalent
(f.e., $\stackrel{f.e.}{\sim}$)} if there
exist matrices $P,Q,K$ of suitable sizes such that
$(A',B')=(P(A+BK)P^{-1},PBQ)$, with $P,Q$ invertible.

Note that equivalence and similarity
of matrices
over $R$ are particular cases of feedback equivalence,
indeed, if $(A,B)$ and $(A',B')$ are feedback equivalent, then $B,B'$
are equivalent, while the equivalence of $(A,0),(A',0)$ implies similarity
of $A,A'$.

Classification and canonical forms under feedback equivalence for
linear systems are classical problems in linear control theory.
While a general solution for arbitrary rings can not be expected
several specific
results have been obtained for special classes of commutative rings. If $R$ is a field,
controllability indices and Brunovski canonical form for reachable systems together with
Kalman's decomposition give a complete description.

For larger classes of rings a variety of partial results have been obtained so
far, e.g. a partial canonical form for reachable systems over principal ideal domains
\cite{pseudo-normal} or a canonical form for weakly reachable single-input systems over
B\'ezout domains \cite{singleinput}.

In \cite{pointwise}, \cite{andresLAA2010}   it was observed, that certain basic control problems have a solution if and
only if the underlying ring is regular (in the sense of von Neumann, see definition below):
\\ -- local and global
feedback equivalence over $R$ are equivalent iff $R$ is regular.
\\ -- cyclizability of the reachable part of systems is always possible iff $R$ is regular.\\
For details see \cite[Theorem 4.1]{pointwise} and \cite[Theorem 4]{andresLAA2010}.\\
At the same time regular rings have nice linear algebraic properties, for example:\\[-3ex]

\noindent -- $R^{n\times n}$ is regular iff $R$ is regular
\cite{goodearl}.\\
-- The tensor criterion for matrix similarity of Byrnes and Gauger
is valid iff the underlying ring is regular
\cite{gustafson}.

This encourages to study feedback classification and canonical forms
over the large class of commutative regular rings.
A commutative ring $R$ is said to be (von Neumann) {\it regular} if for any $a$ in $R$
there exists $x\in R$ such that $a^2x=a$.

The following mostly well-known properties of regular rings will be of use for us later,
see \cite{goodearl} for details.\\[-5ex]

\begin{enumerate} \itemsep -1ex
\item[-] Any finitely generated ideal is principal
(i.e. $R$ is a B\'ezout ring) and generated by an idempotent.

\item[-] Any element is the product of a unit with an idempotent.
\item[-]
$R$ has stable range 1, which for a B\'ezout ring means that whenever $aR+bR=dR$ for
some $a,b,d\in R$ then there exists $c\in R$ such that
$a+bc=ud$, with $u$ a unit.
\item[-] $R$ is an elementary divisor ring, which means that there exists a `Smith Form':
every matrix $B$ is equivalent to a diagonal matrix with diagonal entries
$d_1|d_2| \cdots | d_r$.
In fact, the Smith Form is unique, if the $d_i$'s are idempotents
(see Lemma \ref{lemaformareducida} below).
\item[-] Any finitely generated submodule of $R^n$ has a finitely generated complement.
\end{enumerate}

Note also that any field is a regular ring and more generally any subring of any product
of fields and in particular many rings of functions comprising many finite rings.

\section{Preliminaries}
For an $n$-dimensional system $\Sigma=(A,B)$ over a commutative ring the right image
$\mathcal R_\Sigma$ of $[B,AB,\dots,A^{n-1}B]$ is called the {\it module of reachable states}.
$\Sigma$ is {\it reachable} if $\mathcal R_\Sigma=R^n$. When $R$ is regular, then $\mathcal R_\Sigma$
is always a (finitely genarated) direct summand of $R^n$ with a finitely generated complement.
If $\mathcal R_\Sigma$ has a finite basis which can be completed to a basis of $R^n$,
then the system $\Sigma$ is f.e. to a system of the form
\[
    \big(
    \bmx
        A'&A''\\0&A'''
    \emx,
    \bmx
        B'\\0
    \emx
    \big)
\]
where $(A',B')$ is reachable. This is usually called a {\it Kalman decomposition},
which we shall call {\it strong} in case $A''=0$.
The reachable system $(A,B)$ is in so called {\it Brunovski form} if
\[
A={\sf diag }(A_1,\dots,A_s),B=\big[{\sf diag}(B_1,\dots,B_s),0\big]
\]
where $A_i=\big[e^{(2)},\dots,e^{(n_i)},0\big],B_i=e^{(1)}$
with $e^{(j)}$ the $j$-th canonical basis vector in $R^{n_i}$ for $1\le i\le s$ and
where $n_1\ge \cdots\ge n_s,\ n_1+\cdots +n_s=n$.
The $n_i$ are called controllability or also Kronecker indices. They are a complete
set of invariants under feedback in case $R$ is a field. See e.g. \cite[chapters 3.3 and 4.2 ]{Sontag}

Typically regular rings contain many idempotent elements.
Therefore we now collect some of their  basic properties.
\begin{lema}
\label{isomorfismoproductodirecto}
Let $R$ be a commutative ring and $e_1,\ldots,e_k$ pairwise orthogonal idempotents
($e_ie_j=0$ if $i\neq j$) and such that $\sum_{i=1}^k e_i=1$.
Then:
$$
R =  e_1 R \oplus \cdots \oplus e_k R
$$
Note that the ideal $e_i R$ is at the same time a regular ring with  $e_i$ as unit element.

\end{lema}
%
%

An immediate consequence of such a decomposition is, that
any matricial equation $M=M'$ holds over $R$ if and only if
$e_iM=e_iM'$ holds for all $i$, in particular:
\begin{lema}\label{segmentation}

\begin{itemize}\itemsep -1ex
\item[-] A matrix $P$ is invertible over $R$ iff $e_iP$ is invertible
over $e_iR$ for all $i$.
\item[-] Two matrices $A,A'$ are similar (resp. equivalent) over $R$
iff $e_iA,e_iA'$ are similar (resp. equivalent) over $e_iR$ for all $i$.
\item[-] Two systems $(A,B)$ and $(A',B')$ are f.e. over $R$
iff $(e_iA,e_iB)$ and $(e_iA',e_iB')$ are f.e. over $e_iR$ for all $i$.
\end{itemize}
\end{lema}

Next, we recall basic facts about feedback equivalence.

\begin{lema}
\label{reduccionequivalencia}
Let $R$ be a commutative ring with 1, and
consider the system $(A,B)$ over $R$ of size $(n,m), 0<m< n$ given by
$$
A=
\left[
\begin{matrix}
0   &0 \\
B_1 &A_1\\
\end{matrix}
\right],
B=\left[\begin{matrix} I\\0\end{matrix} \right],
$$
with $I\in R^{m\times m}$ an identity matrix,
$A_1\in R^{(n-m)\times (n-m)}$ and the remaining blocks of
appropriate sizes. Then, one has:
\begin{enumerate} \itemsep -1ex
\renewcommand{\labelenumi}{(\roman{enumi})}
\item
The feedback equivalence class of $(A,B)$ is uniquely determined
by the feedback equivalence class of the $(n-m,m)$ system $(A_1,B_1)$
in the following sense:
two systems
$$
\big(
\left[
\begin{matrix}
0   &0 \\
B_1 &A_1\\
\end{matrix}
\right],
\left[\begin{matrix} I\\0\end{matrix} \right]
\big)
,\:
\big(
\left[
\begin{matrix}
0   &0 \\
B'_1    &A'_1\\
\end{matrix}
\right],
\left[\begin{matrix} I\\0\end{matrix} \right]
\big)
$$
are equivalent if and only if $(A_1,B_1), (A'_1,B'_1)$ are equivalent.
\item
If $R$ is a principal ideal domain or has stable range 1
(in particular, if $R$ is regular),
(i) is also valid if $B$ has any number of additional zero columns,
i.e. it has the block form:
$$
B=\left[\begin{matrix} I & 0 \\0&0 \end{matrix} \right]
$$
\item
If $(A_1,B_1)$ is a reachable system in Brunovski
canonical form then also $(A,B)$ is reachable and can be
transformed
just by  permutation of rows and columns, to Brunovski form.
Moreover the controllability indices
$\kappa_1\ge \cdots \ge \kappa_m,\lambda_1\ge\cdots\ge\lambda_m$
(possibly some of the $\kappa_i$ have to be put to zero)
of $(A_1,B_1)$ and $(A,B)$ are related by $\lambda_i=\kappa_i+1$ for $1\le i\le m$.
\end{enumerate}
\end{lema}

\begin{dem}
(i)
The proof is straightforward and can be adapted from
\cite[Lemmas 2.1 and 2.2]{singleinput}.

(ii) See \cite[Propositions 2.4 and 2.5]{montsehermida}.

(iii) When $(A_1,B_1)$ is in Brunovski form, then for some $s$ we have
\[
    Ae^{(1)}=e^{(m+1)},A^2e^{(1)}=e^{(m+2)},\dots,A^se^{(1)}=e^{(m+s)}.
\]
Similar relations are valid for $e^{(2)},\dots,e^{(m)}$.
Thus an appropriate renumbering of the canonical basis vectors gives a
Brunovski form and the relation for the indices is straightforward.

\end{dem}

\section{Feedback classification and canonical forms}

Let $R$ be a regular ring, and consider a system $\Sigma=(A,B)$ of size
$(n,m)$ over $R$. We will always tacitely assume $B\ne 0$.
In the sequel we will construct recursively
a normal form and a complete set of
invariants for the feedback classification of $\Sigma$. The following reduction step
will be essential:

\begin{lema}
\label{lemaformareducida}
By operations of the feedback group any system $(A,B)$ can be assumed to have
the following form
\begin{equation}
\label{formareducida}
A=\big( a_{ij} \big)_{i,j=1}^n
,\:
B=
\left[
\begin{array}{c|c}
\begin{matrix}
d_1 & & \\ & \ddots & \\ & & d_r
\end{matrix}
& 0 \\ \hline 0 & 0\\
\end{array}
\right]
\end{equation}
with
$d_1 | d_s | \cdots | d_r$,where  all the $d_i$'s are nonzero idempotents
and where for $i=1,\ldots,r$, the $i$-th row of $A$
is orthogonal to $d_i$.
Furthermore, the elements $d_1,\ldots, d_r$ are invariant under feedback equivalence.
If $r=n$ or $r=m$ the corresponding 0-blocks will not occur.
\end{lema}

\begin{dem}
    Since $R$ is an elementary divisor ring and any element of $R$ is a
    product of an idempotent with a unit, the form of $B$ is obtained via matrix
    equivalence.
    Moreover, if we denote by $a_{ij}$ the element in position $(i,j)$ of $A$,
    with a suitable feedback action it can be replaced by
    $a_{ij} -d_ia_{ij}=(1-d_i)a_{ij}$,
    which is orthogonal to $d_i$.

For each $i=1,\ldots,r$, the ideal of $R$
generated by all the $i \times i$ minors of $B$ is invariant under equivalence
and thus invariant under feedback equivalence of the system $(A,B)$.
But this ideal is clearly generated by $d_1\cdots d_i=d_i$,
and an ideal in $R$ cannot be generated by two distinct idempotents.
\end{dem}

We are now ready to solve the feedback classification problem
for systems over a von Neumann regular ring $R$.
A first attempt in this direction could be to reduce the problem
to the classification of systems over all residue fields
(see \cite[Theorem 4.1]{pointwise}):
two systems over $R$ are equivalent if and only if
they are equivalent over the residue field
$R/{\mathfrak m}$, for all maximal ideals ${\mathfrak m}$.
However, there may exist an infinite number of maximal ideals in
the regular ring $R$.
But, as we will prove in the next theorem,
the feedback classification problem of a system over $R$
can be reduced to the classification of a {\em finite}
family of systems which {\em behave like}
systems over fields.


\begin{teorema}
[Canonical forms and invariants]
\label{teoremagordo}
Let $\Sigma=(A,B)$ be a system of size $(n,m)$
over a regular ring $R$. Then:
\begin{enumerate} \itemsep -1ex
\renewcommand{\labelenumi}{(\roman{enumi})}
\item There exists a finite family of idempotents $\{e_k\} $,
pairwise orthogonal and with sum 1, such that for each $k$,
the system $(e_kA,e_kB)$ over the regular ring $e_kR$ is feedback
equivalent to a system in strong Kalman decomposition, with the block form
$$
\big(
\left[
\begin{array}{cc}
\widehat{A_k} &0 \\
0 & \widehat{C_k}
\end{array}
\right]
,
\left[
\begin{array}{c}
\widehat{B_k}\\
0
\end{array}
\right]
\big),
$$
where the pair $(\widehat{A_k},\widehat{B_k})$ is reachable and
in Brunovski canonical form over $e_kR$.
\item
With these notations, a complete set of feedback invariants for
$(A,B)$
consists of:
\begin{enumerate} \itemsep -0.5ex
\item [-] The idempotents $\{e_k\} $, which can be obtained successively from the
invariant factors of a finite number of matrix equivalences.
\item [-] The similarity classes of the matrices $\widehat{C_k}$.
\item [-] The controllability indices of the systems
$(\widehat{A_k},\widehat{B_k})$.
\end{enumerate}
\end{enumerate}
\end{teorema}

\begin{dem}
We will proceed by induction on $n$.

$\bullet$
Proof of (i) and (ii) in the case $n=1$.
We can assume without loss of generality
that the system $\Sigma$ is given in the form (\ref{formareducida}):
$A=[a], B=[d,0 \cdots 0]$,
with $d$ idempotent and $ad=0$.
Let us consider the partition $R=dR \oplus (1-d)R$,
and build the systems:
$$
(1-d)\Sigma = ( [(1-d)a], [0 \cdots 0])
\quad \mbox{and }\;
d\Sigma=( [0], [d, 0\cdots 0])=d([0],[1,0,\dots,0])
$$
It is clear that both systems are decomposed in strong Kalman form,
the first one without reachable part, and the second one
already in Brunovski form over the ring $dR$, in which $d$ is the 1-element.
Moreover, the element $d$ is obtained from the equivalence of $B$
and is invariant by Lemma \ref{lemaformareducida},
making $1-d$ also invariant.
Now, by Lemma \ref{segmentation}, the feedback equivalence class of
$\Sigma$ is completely determined by that of $d\Sigma$
(with a single invariant $d$), and that of $(1-d)\Sigma$,
where the first entry $(1-d)a$ in principle can be reduced by similarity,
but in the $1\times 1$-case this has no effect.

$\bullet$
Proof of (i) for $n>1$.
Let $\Sigma$ be given as in (\ref{formareducida}),
and consider the following family of idempotents
(pairwise orthogonal and with sum 1):
$$
\left\{
\begin{array}{rcl}
e_0 &=& 1-d_1\\
e_1 &=&d_1(1-d_2)\\
& \vdots & \\
e_{r-1} &=& d_{r-1}(1-d_r)\\
e_r &=& d_r
\end{array}
\right.
$$

For $i=0$, we have that $e_0=1-d_1$ is orthogonal to $d_1$, and
hence also orthogonal to all the remaining $d_i$'s, so that
$e_0B=0$, which means that $e_0\Sigma$ trivially is
a strong Kalman decomposition with no reachable part.
This step is omitted if $d_1=1$.

For $i=1,\ldots, \min\{r,n-1\}$, note that $e_i$ is a multiple of $d_1,\ldots,d_i$ and
orthogonal to $d_{i+1},\ldots, d_r$ and to the rows $1,\ldots,i$ of $A$.
Note also, that $e_i=0$ if $d_i=d_{i+1}$.

Therefore we obtain for $1\le i\le \min\{r,n-1\}$ and if $e_i\neq 0$:
$$
e_i\Sigma=(e_iA,e_iB)= \big(\left[
    \begin{array}{c|c}
        0 & 0 \\ \hline
        e_iB_i & e_i A_i
    \end{array}
\right],
    \left[
        \begin{array}{c|c}
              e_iI & 0 \\
              \hline
              0 & 0\\
        \end{array}
    \right]\big)
$$
with $A_i \in R^{(n-i)\times (n-i)}$,
$I$ an $i\times i$ identity block,
and the remaining blocks of appropriate sizes.

Since $(e_iA_i,e_iB_i)$ is of size $(n-i,i)$ (with $n-i<n$),
by the induction assumption one can assume that there
exists a finite partition $e_iR=\bigoplus_j e_{ij}R$,
such that for each $j$, the system $(e_{ij}A_i, e_{ij}B_i)$
over the ring $e_{ij}R$ is equivalent to one of the form
$$
\big(
\left[
\begin{array}{cc}
\widetilde{A_{ij}} &0 \\
0 & \widetilde{C_{ij}}
\end{array}
\right]
,
\left[
\begin{array}{c}
\widetilde{B_{ij}}\\
0
\end{array}
\right]
\big),
$$
with $(\widetilde{A_{ij}},\widetilde{B_{ij}})$ reachable and in Brunovski
canonical form.

In virtue of Lemma \ref{reduccionequivalencia},(i),(ii), we obtain
$$
e_{ij}\Sigma=(e_{ij}A,e_{ij}B) \stackrel{f.e.}{\sim}
\big(
\left[
\begin{array}{c|cc}
0 & 0 &0 \\ \hline\\[-2.2ex]
\widetilde{B_{ij}} &\widetilde{A_{ij}} &0 \\
0 &0 &\widetilde{C_{ij}}
\end{array}
\right]
,\;
\left[
\begin{array}{c|c}
e_{ij}I & 0 \\ \hline
0 &0 \\
0 &0
\end{array}
\right]
\big)
$$
The system on the right side can also be written as
$$\big(
\left[
\begin{array}{cc|c}
0 & 0 &0 \\
\widetilde{B_{ij}} &\widetilde{A_{ij}} &0 \\\hline\\[-2.2ex]
0 &0 &\widetilde{C_{ij}}
\end{array}
\right]
,\;
\left[
\begin{array}{cc}
e_{ij}I & 0 \\
0 &0 \\\hline
0 &0
\end{array}
\right]
\big)\ .
$$

By Lemma \ref{reduccionequivalencia},(iii) we know that

$$
\big(
\left[
\begin{matrix}
0 &0 \\
\widetilde{B_{ij}} &\widetilde{A_{ij}}
\end{matrix}
\right]
,
\left[
\begin{matrix}
e_{ij}I &0\\
0 &0
\end{matrix}
\right] \big) \stackrel{f.e.}{\sim}
\big(
\widehat{A_{ij}} , \widehat{B_{ij}
\\
\big)
},
$$ with $\big(\widehat{A_{ij}},\widehat{B_{ij}}\big)$ in Brunovski form.
If we let $\widehat{C_{ij}}=\widetilde{C_{ij}}$, then we obtain finally
\[
    e_{ij}\Sigma\ \stackrel{f.e.}{\sim}\
    \big( \left[\begin{matrix}
                \widehat{A_{ij}}&0\\0&\widehat{C_{ij}}
          \end{matrix}\right],
          \left[\begin{matrix}
                \widehat{B_{ij}}\\0
          \end{matrix}\right]
    \big)\
\]
with the desired properties for the righthand side system.
This completes the study of the systems $e_i\Sigma$,
for $i=1,\ldots,\min\{r,n-1\}$.

If $r=n$, the system $e_n\Sigma$ is already in Brunovski form,
because
$e_nA$ is zero and $e_nB=[e_nI|0]$, with $I$ an $n\times n$ identity block.

To sum up, we have obtained a finite partition of $R$ with
idempotents $\{e_{ij}\}$ such that for each $i,j$, the system
$e_{ij}\Sigma$ is a strong Kalman decomposition
with reachable part in Brunovski form, as required.

$\bullet$
Proof of (ii) for $n>1$.
With the preceeding notations, the elements $\{e_i\}$ are
constructed from the Smith form of $B$ and
satisfy the conditions of Lemma \ref{isomorfismoproductodirecto},
from which it follows by Lemma \ref{segmentation} that the
classification of $\Sigma$ over $R$ is reduced to the classification
of the systems $e_i\Sigma$ over the regular rings $e_iR$
(at this step we need only consider those indices $i$ for
which $e_i\not= 0$).
A complete set of invariants for $\Sigma$ will be given by
the elements $e_0,\ldots,e_r$,
together with all the invariants obtained recursively.

For $i=0$ we have seen in (i) that $e_0B=0$,
i.e. the feedback class of $(e_0A,e_0B)$ is reduced to the
similarity class of $e_0A$.

For $i=1,\ldots, \min\{r,n-1\}$, by Lemma \ref{reduccionequivalencia}
the feedback class of $e_i\Sigma$ is
completely determined by the feedback class of the system
$(e_iA_i,e_iB_i)$,
which by the induction hypothesis (see the notation in the proof of part (i))
reduces to a collection of idempotents $\{e_{ij}\} $
obtained from equivalences of matrices, the similarity class of all matrices
$\widetilde{C_{ij}}$, and the controllability indices of
all the systems $(\widetilde{A_{ij}},\widetilde{B_{ij}})$,
which are directly related to the controllability indices of the augmented systems
$(\widehat{A_{ij}},\widehat{B_{ij}})$ as was explaineed in
\mbox{Lemma \ref{reduccionequivalencia},(iii).}

Finally, if $r=n$, the additional system $e_n\Sigma$
has the trivial form $(0, [e_nI|0])$ and
does not provide any new invariant.
This completes the proof.
\end{dem}

\section{The reachable and the single input case}
As an immediate consequence, reachable systems can be completely
described without reference to similarity, i.e. only with idempotents
and Brunovski blocks.

\begin{corolario}
[Reachable systems]
With the above notations,
every reachable system over a regular ring
is equivalent to a direct sum of systems
in Brunovski canonical form, and the feedback equivalence class can be
described by a family of idempotents and Brunovski indices.
\end{corolario}

\begin{dem}
With the notations of the previous theorem, reachability
of a given system $(A,B)$ is equivalent to reachability of all
the systems
$$
(e_k A,e_k B) \stackrel{f.e.}{\sim}
\big(
\left[
\begin{array}{cc}
\widehat{A_k} &0 \\
0 & \widehat{C_k}
\end{array}
\right]
,
\left[
\begin{array}{c}
\widehat{B_k}\\
0
\end{array}
\right]
\big),
$$
which is only possible if none of the blocks $\widehat{C_k}$ appear,
i.e. one can remove all references to similarity in the previous theorem.
%
%
\end{dem}

The description via a finite collection of idempotents seems essential and unavoidable for an
explicit exposure of some kind of normal form. Nevertheless, as in the classical field case,
over a regular ring it is possible to classify reachable linear systems in a more global,
but at the same time more abstract way, by the following sequence of
submodules of the module $\mathcal R_\Sigma$ of reachable states.

For an $n$-dimensional system  $\Sigma=(A,B)$ over a regular
ring $R$ let
$$
N_k^{\Sigma} = \mathrm{im}[B|AB| \cdots |A^{k-1}B],
$$
for $k=1,\ldots,n$.

\begin{proposicion}[Reachable systems and the modules \pmb{$N_k^\Sigma$}]
\label{equivalenciamodulosN_k}

Let $\Sigma=(A,B)$ and $\Sigma'=(A',B')$ be two reachable systems of size $(n,m)$
over a regular ring $R$. Then, the following statements
are equivalent:
\begin{enumerate} \itemsep -1ex
\renewcommand{\labelenumi}{(\roman{enumi})}
\item $\Sigma$ and $\Sigma'$ are feedback equivalent.
\item The $R^n$-submodules $N_k^{\Sigma}$ and $N_k^{\Sigma'}$
are isomorphic for all $r=1,\ldots,n$.
\end{enumerate}
\end{proposicion}
\begin{dem}
(i) $\Rightarrow$ (ii) is trivial.\ \
(ii) $\Rightarrow$ (i) can be derived directly from results
in \cite{pointwise}.
Let ${\mathfrak m}$ be an arbitrary maximal ideal of $R$.
If we denote by
$\Sigma(\mathfrak{m})$ and $\Sigma'(\mathfrak{m})$ the extensions of
$\Sigma,\Sigma'$ to the residue field $R/\mathfrak{m}$
(by reduction modulo $\mathfrak{m}$), since
$N_k^{\Sigma} \cong N_k^{\Sigma'}$
we must have
$$
N_k^{\Sigma(\mathfrak{m})}\cong N_k^{\Sigma'(\mathfrak{m})}
$$
for $k=1,\ldots,n$.
But the dimensions of the $R/{\mathfrak m}$-vector spaces
$N_k^{\Sigma(\mathfrak{m})}, N_k^{\Sigma'(\mathfrak{m})}$
characterize completely the equivalence classes of the
reachable systems $\Sigma(\mathfrak{m}),\Sigma'(\mathfrak{m})$
over the field $R/{\mathfrak m}$
(see \cite[Proposition 2.5]{pointwise}), therefore it follows
that $\Sigma(\mathfrak{m})\ \pmb{\stackrel{f.e.}{\sim}}\  \Sigma'(\mathfrak{m})$
over $R/{\mathfrak m}$.
Now, by \cite[Theorem 4.1]{pointwise} we can conclude that
$\Sigma \stackrel{f.e.}{\sim} \Sigma'$ over the regular ring $R$.
\end{dem}

Proposition \ref{equivalenciamodulosN_k} is yet another instance where regular rings behave
``classically''.
Even for principal ideal domains
a corresponding general result is not
possible because of a counterexample in \cite[section 6.2]{BSch}.

As a consequence of Proposition \ref{equivalenciamodulosN_k} we obtain

\begin{corolario}
The feedback equivalence class of a reachable system $\Sigma=(A,B)$
of size $(n,m)$ over a regular ring $R$ is completely
determined by the invariant factors of the matrices
$
[B|AB|\cdots |A^{k-1}B]
,$
for $k=1,\ldots,n$.
\end{corolario}
\begin{dem}
The result is clear, since for each $k=1,\ldots,n$, the columns of the matrix
$[B|AB|\cdots |A^{k-1}B]$ generate the $R$-module $N_k^{\Sigma}$,
whose isomorphism class is determined by the invariant factors of the
generating matrix.
\end{dem}

In the single-input case, the reachable part of a system can be
transformed to a generalized controller canonical form. All notations are as before.
\begin{proposicion}[Single input systems]
.\\[-5ex]
\label{impulsosimple}
\begin{enumerate}
\renewcommand{\labelenumi}{(\roman{enumi})}
    \item
        A canonical form for a system $\Sigma=(A,B)$, $B\ne 0$,
        of size $(n,1)$ over a regular ring $R$ is given by
        $$
        \widetilde{A}=
        \left[
        \begin{matrix}
        0 & 0 &\cdots &0 &0 \\
        d_2 &0 &\ddots &\vdots &\vdots \\
        0 &d_3 &\ddots &\vdots &\vdots \\
        \vdots &\ddots & \ddots &0 &0 \\
        0 &\cdots &0 &d_n &0
        \end{matrix}
        \right]
        + A^{\displaystyle *}
        ,\:
        \widetilde{B}=
        \left[
        \begin{matrix}
        d_1\\
        0\\
        \vdots\\
        \vdots\\
        0
        \end{matrix}
        \right]
        $$
        with idempotents $d_1 |d_2 | \cdots |d_n$,
        possibly  $d_{r+1}=\cdots =d_n=0$ for some $r\ge 1$, and with
        \[
        A^{\displaystyle *}=\widetilde{A}_{1-d_1}
        +
        \left[
        \begin{matrix}
        0 & 0 \\
        0 & \widetilde{A}_{d_1(1-d_2)}
        \end{matrix}
        \right]
        +\cdots +
        \left[
        \begin{matrix}
        0 & 0 \\
        0 & \widetilde{A}_{d_{r-1}(1-d_r)}
        \end{matrix}
        \right],
        \]
        where $\widetilde{A}_{d_i(1-d_{i+1})} \in R^{(n-i)\times (n-i)}$ for $1\le i<r$
        and $d_{i+1}\widetilde{A}_{d_i(1-d_{i+1})}=0$ for $0\le i<r$.

    \item For $k=1,\ldots,n$ one has
        $$
        N_k^{(A,B)}\cong N_k^{(\widetilde{A},\widetilde{B})}=d_1e^{(1)}R+\cdots + d_ke^{(k)}R \cong d_1R \oplus d_2R \oplus \cdots \oplus d_kR
        $$
    \item The elements $d_1,\ldots,d_n$ can be obtained as
        the idempotent invariant factors appearing in the Smith form of the matrix
        $[B,AB,\ldots,A^{n-1}B]$.
\end{enumerate}
\end{proposicion}

\begin{dem}
1.: For a system $\Sigma=(A,B)$ of size $(n,1)$,
    the reduced form (\ref{formareducida}) yields only one idempotent invariant factor
    $d_1$, i.e. two idempotents $e_0=1-d_1$ and $e_1=d_1$ as in the proof of
    Theorem \ref{teoremagordo}.
    The system $e_0\Sigma$ can only be reduced further by similarity and is characterized by
    some corresponding canonical form $\widetilde{A}_{1-d_1}$,
    while
    the class of $d_1\Sigma$ is reduced to that of a system
    $(d_1A_1,d_1B_1)$ of size $(n-1,1)$.
    Continuing this iteration $r$ times, at each step $i$
    we will have a system $(d_iA_i,d_iB_i)$ of size $(n-i,1)$,
    with $d_{i+1}$ the idempotent invariant factor of $d_iB_i$ ($d_i| d_{i+1}$),
    and $\widetilde{A}_{d_i(1-d_{i+1})}$ a canonical form of
    $d_i(1-d_{i+1})A_i$ for similarity over the ring $d_i(1-d_{i+1})R$.
    Multiple use of Lemma \ref{reduccionequivalencia},(i)  leads to the stated result.

2. and 3.: The modules $N_k^\Sigma$ are invariant under feedback and remain
    isomorphic under feedback equivalence. Therefore 2. and 3. are straightforward consequences of 1.
\end{dem}

\section{Implementation}
The procedure given in Theorem \ref{teoremagordo} is constructive,
which gives rise to an explicit algorithm, provided Smith normal
forms are computable over the ring $R$, and a canonical form is known
for the similarity relation.

\begin{algoritmo}[Canonical forms].
\label{formascanonicas}
\par\noindent
INPUT: $(e,A,B)$, whith $(A,B)$ a system and $e$ an idempotent
(initially $e=1$)
\par\noindent
OUTPUT: A list of lists $[e_k, \widehat{A_k}, \widehat{B_k}]$,
where $e_k$ is an idempotent and $(\widehat{A_k}, \widehat{B_k})$
is a canonical form of $(A,B)$ over the ring $e_kR$.
\begin{enumerate} \itemsep -0.7ex
\item
Replace $(A,B)$ by its reduced form (\ref{formareducida}),
with $d_1|d_2|\ldots |d_r$ idempotents.
\item
Initialize L=$[\;]$ (this will be the output).
\item
Define $\{ e_0=e-d_1, e_1=d_1(e-d_2), \ldots,
e_{r-1}=d_{r-1}(e-d_r), e_r=d_r$\}.
\item
For each $e_i\neq 0$, do steps 5--8.
\item
If $e_iB=0$, add $[e_i, \widehat{A_{i}}, 0]$ to $L$,
where $\widehat{A_{i}}$ is a canonical form of $A$ for similarity over $e_iR$.
Proceed with next $e_i$.
\item
If $i=r=n$, add
$[e_i, 0, e_iB]$ to $L$.
\item
Extract blocks $e_iA_i,e_iB_i$ from
$
e_iA=
\left[ \begin{matrix} 0&0\\ e_iB_i &e_iA_i \end{matrix}\right].
$
\item
Recursive call with input $(e_i,e_iA_i,e_iB_i)$
and output a list $\left[ e_{ij}, \widetilde{A_{ij}}, \widetilde{B_{ij}}\right]_j$
\item For each $j$, add to $L$:
$
\left[
e_{ij},
\left[\begin{matrix} 0&0\\ \widetilde{B_{ij}} &\widetilde{A_{ij}} \end{matrix}\right],
\left[\begin{matrix} e_{ij}I&0\\ 0 &0 \end{matrix} \right]
\right]
$
\item If $e_i=0$ for all $i$ ($L$ is still $[\; ]$), output
$\big[ [e,A,B]\big] $,
otherwise output $L$.
\end{enumerate}
\end{algoritmo}

\begin{ejem}
{\em
We have used a symbolic software to compute canonical forms
for rings ${\mathbb Z}/(d{\mathbb Z})$, where $d$ is a squarefree integer.
The following system  over the ring ${\mathbb Z}/(210{\mathbb Z})$
splits into 3 subsystems:
\setlength{\arraycolsep}{2pt}
\begin{small}
$$
\begin{array}{c}
\Sigma=
\big(\left[
\begin{matrix}
 21&158&169&147\\
 138&208&43&135\\
 67&46&190&100\\
 167&36&81&203\\
 \end{matrix}
\right]
,
\left[
\begin{matrix}
 178&152&60&58\\
 90&186&36&120\\
 102&96&30&198\\
 140&40&42&146\\
 \end{matrix}
\right]\big)\\
105\Sigma=
\big(
\left[ \begin{matrix}
 105&0&105&105\\
 0&0&105&105\\
 105&0&0&0\\
 105&0&105&105
 \end{matrix}
\right]
,
\left[ \begin{matrix}
 0&0&0&0\\
 0&0&0&0\\
 0&0&0&0\\
 0&0&0&0\\
\end{matrix}
\right]
\big)
,\:
70\Sigma=
\big(
\left[ \begin{matrix}
 0&0&0&0\\
 70&0&0&0\\
 0&70&0&0\\
 0&0&0&140\\
 \end{matrix}
\right]
,
\left[ \begin{matrix}
 70&0&0&0\\
 0&0&0&0\\
 0&0&0&0\\
 0&0&0&0\\
 \end{matrix}
\right]
\big)
\\
36\Sigma=
\big(
\left[ \begin{matrix}
 0&0&0&0\\
 0&0&0&0\\
 0&0&0&0\\
 36&0&0&0\\
 \end{matrix}
\right]
,
\left[ \begin{matrix}
 36&0&0&0\\
 0&36&0&0\\
 0&0&36&0\\
 0&0&0&0\\
 \end{matrix}
\right]
\big)
\end{array}
$$
\end{small}
We see that the classification of $(A,B)$ reduces to the similarity
classification of a $4\times 4$ block of the first system,
and the $1\times 1$ block 140  appearing in the uncontrollable
part of the second system,
while the third system is reachable and completely described
by the sizes of the idempotent and zero blocks which determine
its
canonical form.
E.g, the Brunovski form of $36\Sigma$ is
\[
    36\cdot \big(\left[\begin{array}{cc|c|c}
            0&0&0&0\\
            1&0&0&0\\
            \hline
            0&0&0&0\\
            \hline
            0&0&0&0\\
    \end{array}\right]
    ,
    \left[\begin{array}{c|c|c|c}
    1&0&0&0\\0&0&0&0\\
    \hline
    0&1&0&0\\
    \hline
    0&0&1&0\\
    \end{array}\right]
    \big)
\]
}
\end{ejem}

\section{Concluding Remarks}
Our results show that over a regular ring
feedback classification of systems still behaves quite similar
to the classical context over a field. But typically one has to handle
a finite collection of systems at the same time. Each of these
``parallel''
systems can be transformed constructively into a strong Kalman decomposition, where the reachable part is in
Brunovski canonical form and where the non reachable part can only be altered further via matrix similarity.
Only if a normal form for similarity of matrices is known, as is the case over fields,
our results give a complete normal form.

The connection of our normal form to the elementary divisors of the reachability submodule,
which we state in Proposition \ref{impulsosimple} has not yet been worked out in the multi-input case.

\end{document}